\documentclass[a4paper,11pt]{article}
\usepackage[T1]{fontenc} 
\usepackage{float} 

\makeatletter
\gdef\@fpheader{ }
\gdef\@journal{ }
\makeatother
\RequirePackage{amsmath}
\RequirePackage{amssymb}
\RequirePackage{epsfig}
\RequirePackage{graphicx}
\RequirePackage[numbers,sort&compress]{natbib}
\RequirePackage{color}
\RequirePackage[colorlinks=true
,urlcolor=blue
,anchorcolor=blue
,citecolor=blue
,filecolor=blue
,linkcolor=blue
,menucolor=blue
,pagecolor=blue
,linktocpage=true
,pdfproducer=medialab
,pdfa=true
]{hyperref}

\newif\ifnotoc\notocfalse
\newif\ifemailadd\emailaddfalse
\newif\iftoccontinuous\toccontinuousfalse
\makeatletter
\def\@subheader{\@empty}
\def\@keywords{\@empty}
\def\@abstract{\@empty}
\def\@xtum{\@empty}
\def\@dedicated{\@empty}
\def\@arxivnumber{\@empty}
\def\@collaboration{\@empty}
\def\@collaborationImg{\@empty}
\def\@proceeding{\@empty}
\def\@preprint{\@empty}

\newcommand{\subheader}[1]{\gdef\@subheader{#1}}
\newcommand{\keywords}[1]{\if!\@keywords!\gdef\@keywords{#1}\else%
\PackageWarningNoLine{\jname}{Keywords already defined.\MessageBreak Ignoring last definition.}\fi}
\renewcommand{\abstract}[1]{\gdef\@abstract{#1}}
\newcommand{\dedicated}[1]{\gdef\@dedicated{#1}}
\newcommand{\arxivnumber}[1]{\gdef\@arxivnumber{#1}}
\newcommand{\proceeding}[1]{\gdef\@proceeding{#1}}
\newcommand{\xtumfont}[1]{\textsc{#1}}
\newcommand{\correctionref}[3]{\gdef\@xtum{\xtumfont{#1} \href{#2}{#3}}}
\newcommand\jname{JHEP}

\newcommand\preprint[1]{\gdef\@preprint{\hfill #1}}

\makeatother

\newcommand\note[2][]{%
\if!#1!%
\stepcounter{footnote}\footnotetext{#2}%
\else%
{\renewcommand\thefootnote{#1}%
\footnotetext{#2}}%
\fi}

\makeatletter
\newtoks\auth@toks
\renewcommand{\author}[2][]{%
  \if!#1!%
    \auth@toks=\expandafter{\the\auth@toks#2\ }%
  \else
    \auth@toks=\expandafter{\the\auth@toks#2$^{#1}$\ }%
  \fi
}
\makeatother
\makeatletter
\newtoks\affil@toks\newif\ifaffil\affilfalse
\newcommand{\affiliation}[2][]{%
\affiltrue
  \if!#1!%
    \affil@toks=\expandafter{\the\affil@toks{\item[]#2}}%
  \else
    \affil@toks=\expandafter{\the\affil@toks{\item[$^{#1}$]#2}}%
  \fi
}
\makeatother
\makeatletter
\newtoks\email@toks\newcounter{email@counter}%
\newcommand{\emailAdd}[1]{%
\emailaddtrue%
\ifnum\theemail@counter>0\email@toks=\expandafter{\the\email@toks, \@email{#1}}%
\else\email@toks=\expandafter{\the\email@toks\@email{#1}}%
\fi\stepcounter{email@counter}}
\newcommand{\@email}[1]{\href{mailto:#1}{\tt #1}}
\makeatother

\makeatletter
\newcommand*\collaboration[1]{\gdef\@collaboration{#1}}
\newcommand*\collaborationImg[2][]{\gdef\@collaborationImg{#2}}
\makeatletter
\newcommand\afterLogoSpace{\smallskip}
\newcommand\afterSubheaderSpace{\vskip3pt plus 2pt minus 1pt}
\newcommand\afterProceedingsSpace{\vskip21pt plus0.4fil minus15pt}
\newcommand\afterTitleSpace{\vskip23pt plus0.06fil minus13pt}
\newcommand\afterRuleSpace{\vskip23pt plus0.06fil minus13pt}
\newcommand\afterCollaborationSpace{\vskip3pt plus 2pt minus 1pt}
\newcommand\afterCollaborationImgSpace{\vskip3pt plus 2pt minus 1pt}
\newcommand\afterAuthorSpace{\vskip5pt plus4pt minus4pt}
\newcommand\afterAffiliationSpace{\vskip3pt plus3pt}
\newcommand\afterEmailSpace{\vskip16pt plus9pt minus10pt\filbreak}
\newcommand\afterXtumSpace{\par\bigskip}
\newcommand\afterAbstractSpace{\vskip16pt plus9pt minus13pt}
\newcommand\afterKeywordsSpace{\vskip16pt plus9pt minus13pt}
\newcommand\afterArxivSpace{\vskip3pt plus0.01fil minus10pt}
\newcommand\afterDedicatedSpace{\vskip0pt plus0.01fil}
\newcommand\afterTocSpace{\bigskip\medskip}
\newcommand\afterTocRuleSpace{\bigskip\bigskip}
\newlength{\affiliationsSep}
\newcommand\beforetochook{\pagestyle{myplain}\pagenumbering{roman}}

\DeclareFixedFont\trfont{OT1}{phv}{b}{sc}{11}

\renewcommand\maketitle{
\pagestyle{empty}
\thispagestyle{titlepage}
\setcounter{page}{0}
\noindent{\small\scshape\@fpheader}\@preprint\par

\afterLogoSpace
\if!\@subheader!\else\noindent{\trfont{\@subheader}}\fi
\afterSubheaderSpace
\if!\@proceeding!\else\noindent{\sc\@proceeding}\fi
\afterProceedingsSpace
{\LARGE\flushleft\sffamily\bfseries\@title\par}
\afterTitleSpace
\hrule height 1.5\p@%
\afterRuleSpace
\if!\@collaboration!\else
{\Large\bfseries\sffamily\raggedright\@collaboration}\par
\afterCollaborationSpace
\fi
\if!\@collaborationImg!\else
{\normalsize\bfseries\sffamily\raggedright\@collaborationImg}\par
\afterCollaborationImgSpace
\fi
{\bfseries\raggedright\sffamily\the\auth@toks\par}
\afterAuthorSpace
\ifaffil\begin{list}{}{%
\setlength{\leftmargin}{0.28cm}%
\setlength{\labelsep}{0pt}%
\setlength{\itemsep}{\affiliationsSep}%
\setlength{\topsep}{-\parskip}}
\itshape\small%
\the\affil@toks
\end{list}\fi
\afterAffiliationSpace
\ifemailadd 
\noindent\hspace{0.28cm}\begin{minipage}[l]{.9\textwidth}
\begin{flushleft}
\textit{E-mail:} \the\email@toks
\end{flushleft}
\end{minipage}
\else 
\PackageWarningNoLine{\jname}{E-mails are missing.\MessageBreak Plese use \protect\emailAdd\space macro to provide e-mails.}
\fi
\afterEmailSpace
\if!\@xtum!\else\noindent{\@xtum}\afterXtumSpace\fi
\if!\@abstract!\else\noindent{\renewcommand\baselinestretch{.9}\textsc{Abstract:}}\ \@abstract\afterAbstractSpace\fi
\if!\@keywords!\else\noindent{\textsc{Keywords:}} \@keywords\afterKeywordsSpace\fi
\if!\@arxivnumber!\else\noindent{\textsc{ArXiv ePrint:}} \href{http://arxiv.org/abs/\@arxivnumber}{\@arxivnumber}\afterArxivSpace\fi
\if!\@dedicated!\else\vbox{\small\it\raggedleft\@dedicated}\afterDedicatedSpace\fi
\ifnotoc\else
\iftoccontinuous\else\newpage\fi
\beforetochook\hrule
\tableofcontents
\afterTocSpace
\hrule
\afterTocRuleSpace
\fi
\setcounter{footnote}{0}
\pagestyle{myplain}\pagenumbering{arabic}
} 

\renewcommand{\baselinestretch}{1.1}\normalsize
\divide\@tempdima\baselineskip
\@tempcnta=\@tempdima
\skip\@mpfootins = \skip\footins
\renewcommand{\@dotsep}{10000}

\newcommand\ps@myplain{
\pagenumbering{arabic}
\renewcommand\@oddfoot{\hfill-- \thepage\ --\hfill}
\renewcommand\@oddhead{}}
\let\ps@plain=\ps@myplain

\newcommand\ps@titlepage{\renewcommand\@oddfoot{}\renewcommand\@oddhead{}}

\numberwithin{equation}{section}

\renewcommand\section{\@startsection{section}{1}{\z@}%
                                   {-3.5ex \@plus -1.3ex \@minus -.7ex}%
                                   {2.3ex \@plus.4ex \@minus .4ex}%
                                   {\normalfont\large\bfseries}}
\renewcommand\subsection{\@startsection{subsection}{2}{\z@}%
                                   {-2.3ex\@plus -1ex \@minus -.5ex}%
                                   {1.2ex \@plus .3ex \@minus .3ex}%
                                   {\normalfont\normalsize\bfseries}}
\renewcommand\subsubsection{\@startsection{subsubsection}{3}{\z@}%
                                   {-2.3ex\@plus -1ex \@minus -.5ex}%
                                   {1ex \@plus .2ex \@minus .2ex}%
                                   {\normalfont\normalsize\bfseries}}
\renewcommand\paragraph{\@startsection{paragraph}{4}{\z@}%
                                   {1.75ex \@plus1ex \@minus.2ex}%
                                   {-1em}%
                                   {\normalfont\normalsize\bfseries}}
\renewcommand\subparagraph{\@startsection{subparagraph}{5}{\parindent}%
                                   {1.75ex \@plus1ex \@minus .2ex}%
                                   {-1em}%
                                   {\normalfont\normalsize\bfseries}}

\def\fnum@figure{\textbf{\figurename\nobreakspace\thefigure}}
\def\fnum@table{\textbf{\tablename\nobreakspace\thetable}}

\long\def\@makecaption#1#2{%
  \vskip\abovecaptionskip
  \sbox\@tempboxa{\small #1. #2}%
  \ifdim \wd\@tempboxa >\hsize
    \small #1. #2\par
  \else
    \global \@minipagefalse
    \hb@xt@\hsize{\hfil\box\@tempboxa\hfil}%
  \fi
  \vskip\belowcaptionskip}

\renewenvironment{thebibliography}[1]{%
\begin{oldthebibliography}{#1}%
\small%
\raggedright%
\setlength{\itemsep}{5pt plus 0.2ex minus 0.05ex}%
}%
{%
\end{oldthebibliography}%
}

\begin{document}


\title{\boldmath The Pade Approximant Based Network for Variational Problems
 }


\author[a,b,1]{Chi-Chun Zhou,}\note{zhouchichun@dali.edu.cn}
\author[b,2]{and Yi Liu,}\note{liuyi@dali.edu.cn.}


\affiliation[a]{School of Engineering, Dali University, Dali, Yunnan 671003, PR China}
\affiliation[b]{Department of Physics, Tianjin University, Tianjin 300350, PR China}










\abstract{In solving the variational problem, the key is to
efficiently find the target function that minimizes or maximizes the specified
functional. In this paper, by using the Pade approximant, we suggest a methods
for the variational problem. By comparing the method with
those based on the radial basis
function networks (RBF), the multilayer perception networks (MLP), and the
Legendre polynomials, we show that the method searches
the target function effectively and efficiently.}


\maketitle
\flushbottom


\section{Introduction}

Many problems arising as variational problems, such as the principle of
minimum action in theoretical physics and the optimal control problem in
engineering. In solving the variational problem, the key is to efficiently
find the target function that minimizes or maximizes the specified functional.

In using non-analytical methods to search for the target function,
one faces two problems:
firstly, to ensure that the target function 
is in the range of searching. 
Secondly, to ensure that the target function 
can be found under limited computing power and time.
The first one is a problem of effectiveness and the second is a problem of efficiency

The direct method of Ritz and Galerkin
\cite{gelfand2000calculus,elsgolc2012calculus,giaquinta2013calculus} tries to
express the target function as linear combinations of basis functions and 
reduces the problem to that of solving equations of coefficients. However,
the basis function is determined by the boundary condition, as a result,
the target function might not be expressed by the basis function
and thus be excluded in the range of searching.

One, to improve the method along this line, 
uses Walsh functions \cite{chen1975walsh}, 
orthogonal polynomials \cite{chang1983shifted,horng1985shifted,hwang1983laguerre,razzaghi2000legendre}, 
and fourier series \cite{razzaghi1988fourier,hsiao2004haar} that are
complete and orthogonal to express the target function 
and converts the boundary condition into a constraint of the coefficient. 
In that approach, it is the completeness of the basis function
that ensures the effectiveness. 

Recently, the multilayer perception networks (MLP) 
is used to solve the variational problem
\cite{weinan2018deep,lopez2008neural,gonzalez2009neural}. 
At this case, the functional becomes the loss function and the boundary
condition usually becomes an extra term added to the loss function. 
The MLP is trained to learn the shape of the target function and meet
the boundary condition simultaneously.
In that approach, it is the universal
approximation property of the MLP
\cite{hornik1991approximation,white1990connectionist,hornik1990universal,leshno1993multilayer,hornik1989multilayer}
that guarantees the effectiveness. 

In a word, the completeness of the basis function or the universal
approximation property of the neural networks already solves the problem
of effectiveness. Now the problem of efficiency is to be considered. For example,
in searching the target function with MLP, 
the network might fall into the local minimum instead of the global
minimum.

In this paper, by using the Pade approximant, we suggest a methods
for the variational problem. By comparing the method with 
those based on the radial basis
function networks (RBF), the multilayer perception networks (MLP), and the
Legendre polynomials,  we show that the method searches
the target function effectively and efficiently.

This paper is organized as following. In Sec. 2, we introduce the main method,
where the effective expression of the target function is constructed. In Sec.
3, we solve the illustrative examples. Conclusions and outlooks are given in Sec. 4.

\section{The main method}

In this section, we show the detail of constructing an efficient expresion of
the target function based on the Pade approximant.

\subsection{the Pade approximant: a brief review}

The Pade approximant is a rational function of
numerator degree $m$ and denominator degree $n$
\cite{baker1996pade,brezinski2012history},%

\begin{equation}
y_{pade}\left(  x\right)  =\frac{\sum_{j=1}^{m}w_{j}x^{j}+b_{1}}{\sum
_{i=1}^{n}w_{i}^{\prime}x^{i}+b_{2}},\label{1}%
\end{equation}
where $w_{j}$, $w_{i}^{\prime}$, $b_{1}$, and $b_{2}$ are parameters. For the
sake of convenience, we denote the structure of the Pade approximant as
Pade-$\left[  m/n\right]  $.

Normally, the Pade-$\left[  m/n\right]  $ ought to fit a power series through
the orders $1$,$x$,$x^{2}$,$\ldots$,$x^{m+m}$ \cite{baker1996pade}, that is
\begin{equation}
\sum_{i=0}^{\infty}c_{i}x^{i}=\frac{\sum_{j=1}^{m}w_{j}x^{j}+b_{1}}{\sum
_{i=1}^{n}w_{i}^{\prime}x^{i}+b_{2}}+O\left(  x^{m+n+1}\right)  .
\end{equation}

It is an highly efficient tool to approximate a complex real function with
only countable parameters and has wide applications in many problems
\cite{brent1980fast,cochelin1994asymptotic,langhoff1970pade,loeffel1969pade,vidberg1977solving}%
. Here, we use it as an approximator for the target function.

\subsection{The RBF, the MLP, and the Legendre Polynomial: brief reviews}

In order to compare the method with those based on the radial basis function
networks (RBF), the multilayer perception networks (MLP), and the Legendre
polynomials, we give a brief review on the RBF, the MLP, and the Legendre Polynomial.

\textit{The MLP. }The MLP or multilayer feed-forward networks is a typical
feed-forward neural network. It transforms an $n$-dimensional input $x$ to a
$k$-dimensional output $y$ and implements a class of mappings from $R^{n}$ to
$R^{k}$
\cite{hornik1991approximation,white1990connectionist,hornik1990universal,leshno1993multilayer,hornik1989multilayer}%
. The building block of the MLP is neurons where a linear and non-linear
transforms are successively applied on the input. A collection of neurons forms
a layer and a collection of layers gives a MLP. For example, for a one-layer
MLP with the number of hidden node, the neuron in the layer, $l$, the relation
between the input $x$ and the output $y$ can be explicitly written as%

\begin{equation}
y_{mlp}\left(  x\right)  =\sum_{i=1}^{l}w_{i}^{\prime}\sigma\left(  \sum
_{j=1}^{n}w_{ij}x_{j}+b_{1}\right)  +b_{2},\label{2}%
\end{equation}
where $y$ is a $1$-dimensional output in this case, $\sigma\left(  x\right)  $
is the non-linear map called the active function and is usually chosen to be
$\tanh x$ or
\begin{equation}
sigmoid\left(  x\right)  =\frac{1}{1+e^{-x}}.
\end{equation}
$w_{ij}$, $w_{i}^{\prime}$, $b_{1}$, and $b_{2}$ are parameters. By tuning the
parameters, the MLP is capable of approximating a target function. For the
sake of convenience, we denote the structure of the MLP as MLP-$[[l,\sigma
\left(  x\right)  ]]$. E.g., a two-layer MLP with $32$ neurons in each layer
and activate functions both $sigmoid\left(  x\right)  $ is 
MPL-$[[32,sigmoid],[32,sigmoid]]$.

\textit{The RBF. }Beyond the MLP,\textit{ }RBF or the radial basis function
networks is another typical neural network. Similarly, it also transforms an
$n$-dimensional input $x$ to a $k$-dimensional output $y$ and implements a
class of mappings from $R^{n}$ to $R^{k}$
\cite{park1991universal,orr1996introduction,park1993approximation}. However,
the structure of the RBF is different: the distance between the input and the
center are transformed by a kernel function. The linear combination of the
result of the kernel function gives the output $y$. For example, for a RBF
with $l$ centers, the relation between the input $x$ and the output $y$ can be
explicitly written as
\begin{equation}
y_{rbf}\left(  x\right)  =\sum_{j=1}^{l}w_{j}\exp\left[  -\frac{\sum_{i=1}%
^{n}\left(  x_{i}-c_{ji}\right)  ^{2}}{2\sigma_{j}}\right]  +b,\label{3}%
\end{equation}
where, the kernel function is the Gauss function in this case. The RBF is also
a good approximator
\cite{park1991universal,orr1996introduction,park1993approximation} and, at
some cases, more efficient than the MLP. For the sake of convenience, we denote
the structure of the RBF as RBF-$[l]$.

\textit{The Legendre polynomial. }In real analysis, a real function can be
expressed as a linear combination of basis such as complete polynomials
\cite{dunkl2014orthogonal}. The Legendre polynomial is complete and
orthogonal. It satisfies the recurrent relations \cite{weisstein2002legendre}
\begin{equation}
P_{n+1}\left(  x\right)  =\frac{2n+1}{n+1}xP_{n}\left(  x\right)  -\frac
{n}{n+1}P_{n-1}\left(  x\right)
\end{equation}
for $n=1$, $2$, $3$, $\ldots$, where $P_{j}\left(  x\right)  $ is Legendre
polynomial of order $j$, $P_{0}\left(  x\right)  =1$, and $P_{1}\left(
x\right)  =x$. In this work, we express the target function as
\begin{equation}
y_{legend}\left(  x\right)  =\sum_{j=1}^{m}w_{j}P_{j}\left(  x\right)
+b\label{5}%
\end{equation}
with $w_{j}$ and $b$ being parameters. For the sake of convenience, we denote
the structure as Leg-$m$.

\textit{The power polynomial. }For the reader, the power polynomial is a
familiar tool to approximate a function. For example, the Taylor expansion, a
textbook content, is based on it. Here, in order to show that the methods such
as the MLP and the RBF are nothing mysterious but merely an approximator, we
give the result based on the power polynomial. We express the target function
as
\begin{equation}
y_{poly}\left(  x\right)  =\sum_{j=1}^{m}w_{j}x^{j}+b,\label{4}%
\end{equation}
where $w_{j}$ and $b$ are parameters. We show that the neural-network method
differs from the power-polynomial method only in efficiency. For the sake
of convenience, we denote the structure as Poly-$m$.

\subsection{The expression of the target function}
In searching the target function that minimizes or maximizes the
specified functional numerically, the parameter is tuned to shape the output function.
However, the boundary condition might reduce the efficiency, because it
becomes an extra constraint on the parameter, i.e., the parameter now is tuned
not only to shape the function, but also to move the function to the fixed
point. In this section, we suggest a expression for the target function
which has the universal approximation property and satisfies the boundary condition
automatically. With this approach, the boundary condition is no more an extra
constraint on the parameter. 

There are various kinds of boundary conditions in the variational problem.
Here, without loss of generality, we focus on  $1$-dimensional problems with
the fix-end boundary condition.

\textit{The boundary factor}. We introduce the boundary factor,
\begin{equation}
bound\left(  x\right)  =\left(  x-x_{a}\right)  ^{m_{a}}\left(  x_{b}%
-x\right)  ^{m_{b}},\label{21}%
\end{equation}
where $m_{a}$ and $m_{b}$ are parameters. $x_{a}$ and $x_{b}$ are boundaries
of $x$. We, for the sake of convenience, rewrite the output of Eqs. (\ref{1})-(\ref{4})
as
\begin{equation}
y_{net}\left(  x\right)  =\left\{
\begin{array}
[c]{c}%
y_{pade}\left(  x\right)  \\
y_{MLP}\left(  x\right)  \\
y_{RBF}\left(  x\right)  \\
y_{legendre}\left(  x\right)  \\
y_{poly}\left(  x\right)
\end{array}
\right.  .\label{22}%
\end{equation}
Multiplying the boundary factor,
Eq. (\ref{21}), to the output $y_{net}\left(  x\right)  $ ensures 
that the output function passes through points $\left(
x_{a},0\right)  $ and $\left(  x_{b},0\right)  $.

\textit{The construction. }In order to pass through the fix-end points $\left(
x_{a},y_{a}\right)  $ and $\left(  x_{b},y_{b}\right)  $, we add a function
\begin{equation}
g\left(  x\right)  =x\ast\frac{y_{b}-y_{a}}{x_{b}-x_{a}}+\frac{x_{b}%
y_{a}-y_{b}x_{a}}{x_{b}-x_{a}}%
\end{equation}
to the output. Finally, the expression of the target function reads%
\begin{equation}
y_{net\_final}\left(  x\right)  =y_{net}\left(  x\right)  \ast bound\left(
x\right)  +g\left(  x\right)  .\label{23}%
\end{equation}
Eq. (\ref{23}) inherits the good approximate ability from the Pade
approximant, the MLP and so on and passes through the fixed-end point simultaneously.

\subsection{The loss function and the learn algorithm}

\textit{The loss function. }The specified functional now read%
\begin{equation}
J=\int_{x_{a},y_{a}}^{x_{b},y_{b}}dxF\left[  y_{net\_final}\left(  x\right)
,y_{net\_final}^{\prime}\left(  x\right)  ,\ldots\right]  ,\label{31}%
\end{equation}
where $F$ is the specified function of $y_{net\_final}\left(  x\right)  $,
$y_{net\_final}^{\prime}\left(  x\right)  $, and so on. In order to conduct a
numerical computation, Eq. (\ref{31}) is approximated by a summation
\begin{equation}
loss=\frac{\left(  x_{b}-x_{a}\right)  }{N}\sum_{i=1}^{N}F\left[
y_{net\_final}\left(  x_{i}\right)  ,y_{net\_final}^{\prime}\left(
x_{i}\right)  ,\ldots\right]  ,\label{32}%
\end{equation}
where $N$ is the number of sample points, $x_{i}$ is sampled uniformly from
$\left(  x_{a},x_{b}\right)  $. Thus the variational problem converted into a
optimization problem:%
\begin{equation}
\min_{w,b,m}loss.
\end{equation}

\textit{The gradient descent method and the back-propagation algorithm. }We
use the gradient descent method to find the optimal parameter, e.g., the
parameter is updated by the following equation
\begin{equation}
w_{new}=w_{old}+l\left.  \frac{\partial loss}{w}\right\vert _{w=w_{old}%
},\label{33}%
\end{equation}
where $l$ is the learning rete. $w_{old}$ and $w_{new}$ are the old and new
parameters after one step respectively. In Eq. (\ref{33}), $\frac{\partial
loss}{w}$ is calculated by the back-propagation algorithm%
\begin{equation}
\frac{\partial loss}{w}=\frac{\partial loss}{\partial y_{net\_final}}%
\frac{\partial y_{net\_final}}{\partial w}+\frac{\partial loss}{\partial
y_{net\_final}^{\prime}}\frac{\partial y_{net\_final}^{\prime}}{\partial
w}+\ldots.
\end{equation}

An implementation based on python and tensorflow is given in 
\href{https://github.com/zhouchichun/mlp_rbf_poly_legend_pade}{github}. In the
implementation, the back-propagation algorithm is automatically processed and
the Adam algorithm, a developed gradient descent method is applied. 

\section{The illustrative example}

In this section, we use the method to solve variational problems that are
partly collected from the literatures
\cite{chang1983shifted,horng1985shifted,hwang1983laguerre,razzaghi2000legendre,razzaghi1988fourier,hsiao2004haar,weinan2018deep,lopez2008neural,gonzalez2009neural}%

1) \textit{The shortest path problem. }The functional reads\textit{ }
\begin{equation}
J=\int_{-1}^{1}\sqrt{1+\left(  y^{\prime}\right)  ^{2}}dx\label{e1}%
\end{equation}
with boundary condition%
\begin{equation}
y\left(  -1\right)  =0\text{ and }y\left(  1\right)  =2.\label{e1b}%
\end{equation}
Exact results are%
\begin{align}
y_{exact}\left(  x\right)    & =x+1,\nonumber\\
J\left(  y_{exact}\right)    & =2\sqrt{2}\simeq2.8284.
\end{align}
The $y_{exact}\left(  x\right)  $ is a straight line at this case, 
however, it dose not mean that the
task is simple, because to
find the target function without any pre-knowledge is much more difficult than
to learn to express a known target function.

Numerical results are%

\begin{figure}[H]
\centering
\includegraphics[width=0.8\textwidth]{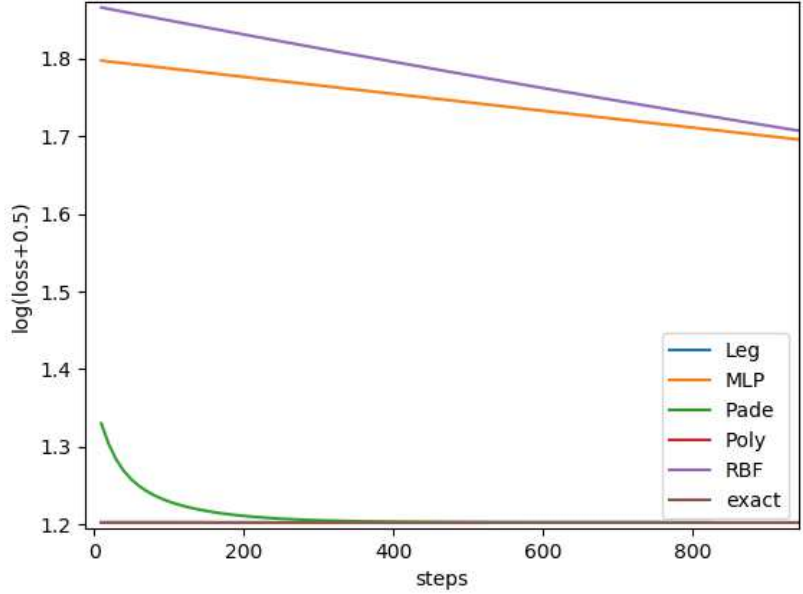}
\caption{The loss versus the steps of different methods.}
\label{shortest}
\end{figure}

\begin{table}[H]  
 \label{table}
\centering
 \begin{tabular}{llll}  

\hline   

  $\text{Structure }$ & $number of parameters$ & $J\left(  y_{net\_final}\right)$ &$\text{Relative error}$\\  

\hline   
  $\text{Pade-}[5/5]$ & $12$& $2.8285$ & $-3.5\times10^{-5}$  \\ 
  $\text{RBF-}[8]$& $25$& $2.8285$ & $-3.5\times10^{-5}$  \\
  $\text{MLP-}[[8,sigmoid]] $& $18$& $2.8285$ & $-3.5\times10^{-5}$  \\
  $\text{Leg-}10$& $11$& $2.8274$ & $-3.5\times10^{-4}$ \\
  $\text{Poly-}10 $& $11$& $2.8285$ & $-3.5\times10^{-5}$  \\

\hline 
\end{tabular}
\end{table}
The efficiency of each methods are show in Fig. (\ref{shortest})

From Fig. (\ref{shortest}), one can see that the method based on the Pade approximant
converges faster than those based on the RBF and the MLP.
In the method based on Legendre polynomials and the power polynomials
the initial value happens to be the target function. 

2) \textit{The minimum drag problem. }The functional reads%
\[
J=\int_{0}^{1}yy^{\prime3}dx
\]
with boundary condition%
\begin{equation}
y\left(  0\right)  =0\text{ and }y\left(  1\right)  =1.
\end{equation}

Exact results are%

\begin{align}
y_{exact}\left(  x\right)    & =x^{3/4},\nonumber\\
J\left(  y_{exact}\right)    & =\frac{27}{64}\simeq0.4219.
\end{align}

\begin{figure}[H]
\centering
\includegraphics[width=0.8\textwidth]{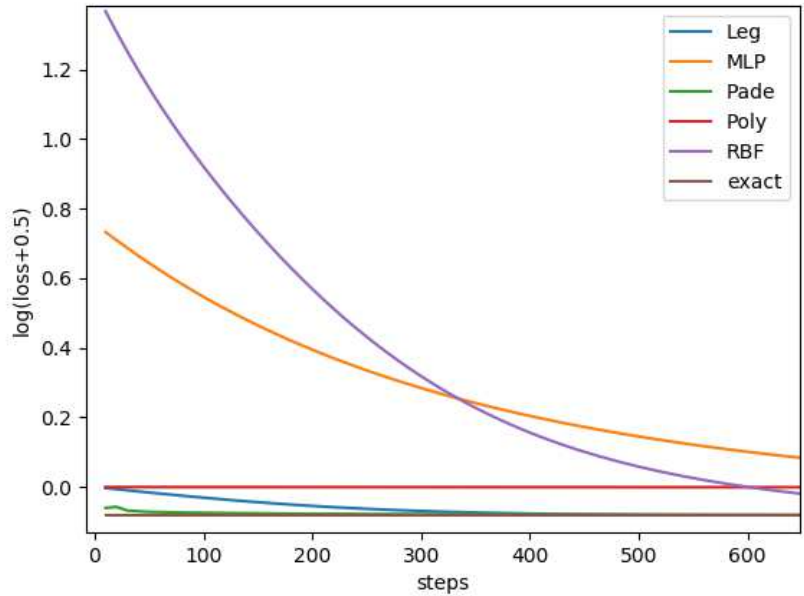}
\caption{The loss versus the step of different methods.}
\label{drag}
\end{figure}

Numerical results are%
\begin{table}[H]
\label{table1}
\centering
 \begin{tabular}{llll}  

\hline   

  $\text{Structure }$& $number of parameters$ & $J\left(  y_{net\_final}\right)$ &$\text{Relative error}$\\  

\hline   
  $\text{Pade-}8/10$& $20$ & $0.4216$ & $7\times10^{-4}$\\
$\text{RBF-}[8] $& $25$& $0.4217$ &  $4.7\times10^{-4}$\\
$\text{MLP-}[[16,sigmoid]]$& $34$ & $0.4220$ & $-2.3\times10^{-4}$\\
$\text{Leg-}15$& $16$ & $0.4218$ & $2.3\times10^{-4}$\\
$\text{Poly-}15$& $16$ & $0.4216$ & $7\times10^{-4}$\\

\hline 
\end{tabular}
\end{table}
The efficiency of each methods are show in Fig. (\ref{drag})
From Fig. (\ref{drag}), one can see the method based on the Pade approximant
again converges faster. The method based on the power polynomials
converges very slow at this time. 
Since we have shown that the method is capable to find the target function,
in the later examples, we eliminate this method because it converges so slow.

3) \textit{A popular illustrative example.} The functional reads%
\[
J=\int_{0}^{1}\left(  y\prime^{2}+xy^{\prime}\right)  dx
\]
with boundary condition%
\begin{equation}
y\left(  0\right)  =0\text{ and }y\left(  1\right)  =\frac{1}{4}.
\end{equation}

Exact results are%

\begin{align*}
y_{exact}\left(  x\right)    & =\frac{1}{2}x\left(  1-\frac{1}{2}x\right)  ,\\
J\left(  y_{exact}\right)    & =\frac{5}{3}\simeq1.667
\end{align*}

\begin{figure}[H]
\centering
\includegraphics[width=0.8\textwidth]{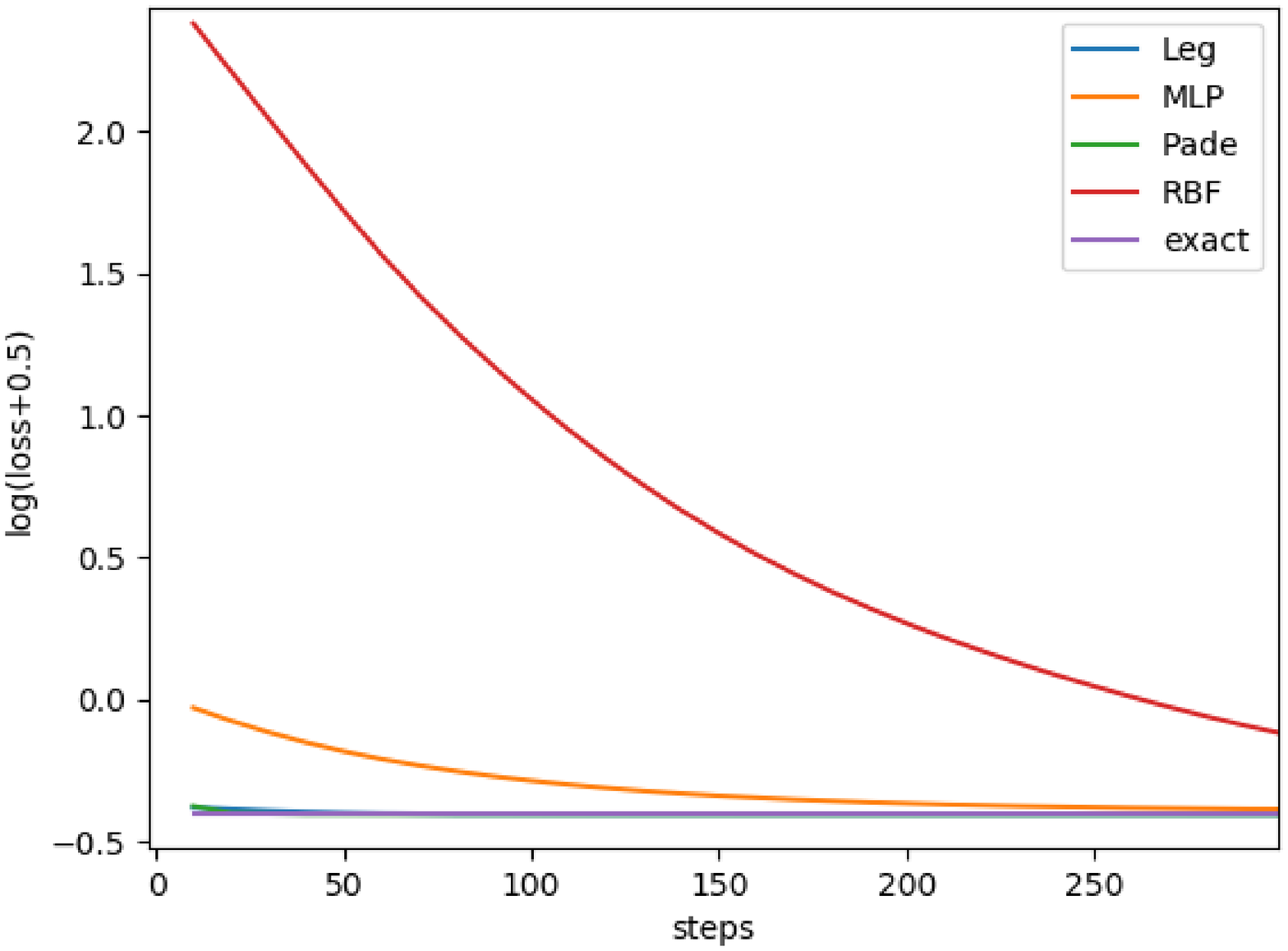}
\caption{The loss versus the steps of different methods.}
\label{pop}
\end{figure}

Numerical results are%
\begin{table}[H]
\label{table2}
\centering
 \begin{tabular}{llll}  

\hline   

  $\text{Structure }$& $number of parameters$ & $J\left(  y_{net\_final}\right)$ &$\text{Relative error}$\\  

\hline   
$\text{Pade:}8/10 $& $20$& $1.665$ & $1\times10^{-3}$\\
$\text{RBF:}[16] $& $49$& $1.667$ & $0$\\
$\text{MLP:}[[16,sigmoid]] $& $34$& $1.665$&$1\times10^{-3}$\\
$\text{Leg:}15$ & $16$& $1.664$ & $2\times10^{-3}$\\

\hline 
\end{tabular}
\end{table}
The efficiency of each methods are show in Fig. (\ref{pop})

4)\textit{Example 4.}The functional reads%
\[
J=\int_{-\pi/2}^{\pi/2}\left[  y^{\prime2}-2y\cos\left(  x+\frac{\pi}%
{2}\right)  \right]  dx
\]
with boundary condition%
\begin{equation}
y\left(  \pi/2\right)  =0\text{ and }y\left(  -\pi/2\right)  =0.
\end{equation}

Exact results are%
\begin{align}
y_{exact}\left(  x\right)    & =\cos\left(  x+\frac{\pi}{2}\right)  +\frac
{2}{\pi}x,\nonumber\\
J\left(  y_{exact}\right)    & =\frac{4}{\pi}-\frac{\pi}{2}\simeq-0.2976.
\end{align}

\begin{figure}[H]
\centering
\includegraphics[width=0.8\textwidth]{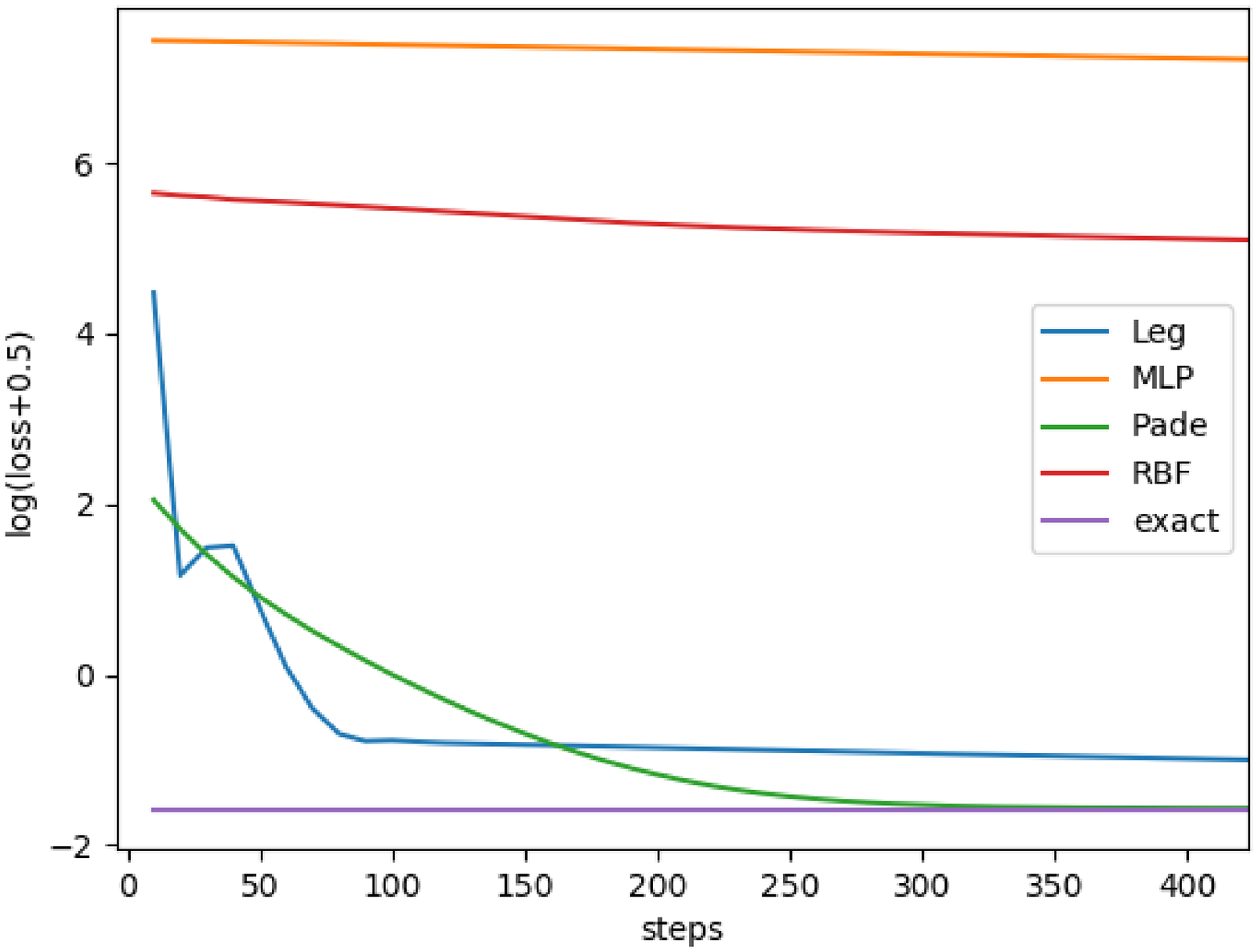}
\caption{ The loss versus the steps of different methods.}
\label{ex4}
\end{figure}

Numerical results are%
\begin{table}[H]
\label{table3}
\centering
 \begin{tabular}{llll}  

\hline   

  $\text{Structure }$& $number of parameters$ & $J\left(  y_{net\_final}\right)$ &$\text{Relative error}$\\  

\hline   
$\text{Pade:}4/5$ & $11$& $-0.2969$ & $-2\times10^{-3}$\\
$\text{RBF:}[16]$& $49$ & $-0.2968 $& $-3\times10^{-3}$\\
$\text{MLP:}[[16,sigmoid]]$& $34$ & $-0.2968$ & $-3\times10^{-3}$\\
$\text{Leg:}15$& $16$ & $-0.2965$ & $-4\times10^{-3}$\\

\hline 
\end{tabular}
\end{table}
The efficiency of each methods are show in Fig. (\ref{ex4}) 

5)\textit{Example 5.} The functional reads%
\[
J=\int_{0}^{1}\left(  y^{\prime2}-y^{2}-2xy\right)  dx
\]
with boundary condition%
\begin{equation}
y\left(  0\right)  =0\text{ and }y\left(  1\right)  =0.
\end{equation}
Exact results are%
\begin{align}
y_{exact}\left(  x\right)    & =\frac{\sin x}{\sin1}-x,\nonumber\\
J\left(  y_{exact}\right)    & =\cot1-\frac{2}{3}\simeq-0.0246.
\end{align}

\begin{figure}[H]
\centering
\includegraphics[width=0.8\textwidth]{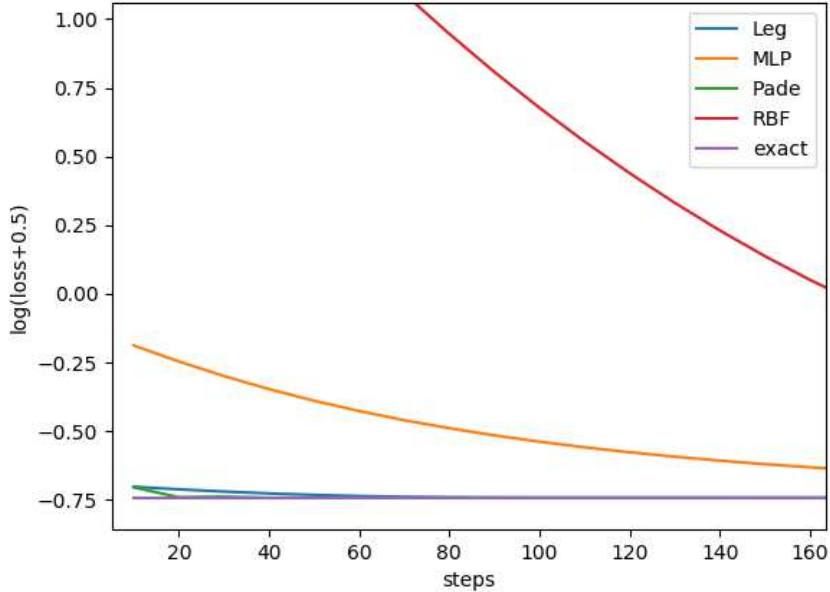}
\caption{The loss versus the steps of different methods.}
\label{ex5}
\end{figure}

Numerical results are%
\begin{table}[H]
\label{table4}
\centering
 \begin{tabular}{llll}  

\hline   
203.
  $\text{Structure }$& $number of parameters$ & $J\left(  y_{net\_final}\right)$ &$\text{Relative error}$\\  

\hline   
$\text{Pade:}4/5 $& $11$&$ -0.0245$ & $-4\times10^{-3}$\\
$\text{RBF:}[16] $& $49$&$ -0.0245 $& $-4\times10^{-3}$\\
$\text{MLP:}[[16,sigmoid]] $& $34$&$ -0.0245$ & $-4\times10^{-3}$\\
$\text{Leg:}15$& $16 $&$ -0.0245$ & $-4\times10^{-3}$\\

\hline 
\end{tabular}
\end{table}
The efficiency of each methods are show in Fig. (\ref{ex5})

\section{Conclusions and outlooks}
In solving the variational problem, the key is to
efficiently find the target function that minimizes or maximizes the specified
functional. Problems of effectiveness and efficiency are both important.
In this paper, by using the Pade approximant, we suggest a methods
for the variational problem. In this approach, the fix-end boundary
condition is satisfied simultaneously.
By comparing the method with
those based on the radial basis
function networks (RBF), the multilayer perception networks (MLP), and the
Legendre polynomials, we show that the method searches
the target function effectively and efficiently.

The method shows that the Pade approximant
can improve the efficiency of neural network.
In solving a many-body system numerically in physics, the efficiency
of the method is important, because the degree of freedom in such system is 
large. The method could be used in searching the wave function
of a many-body system efficiently.
Moreover, it could be used in other tasks, such
as the task of classification and translation.

\section{Acknowledgments}
We are very indebted to Dr Dai for his enlightenment and encouragement.











\bibliographystyle{JHEP}
\bibliography{refs}

\end{document}